\documentclass[12pt,letter]{article}
\usepackage{amssymb,amsbsy,amsmath,amsfonts,amssymb}

\title{Incomparable, non isomorphic and minimal Banach spaces}
\author {Christian Rosendal}
\date {March 2004}
\newcommand {\ca} {{2^\N}}
\newcommand{\cantor}{2^{<\N}}

\newcommand{\ram}{[\N]^{\N}}

\newcommand {\N}{\mathbb N}
\newcommand {\Q}{\mathbb Q}
\newcommand {\R}{\mathbb R}

\newcommand{\compl}{\complement}

\newcommand{\iso}{\cong}

\newcommand{\emb}{\sqsubseteq}

\newcommand{\tom} {\emptyset}

\newcommand{\begr}{\!\upharpoonright}

\newcommand{\bij}{\longleftrightarrow}

\newcommand{\hviss}{\leftrightarrow}

\newcommand{\equi}{\Longleftrightarrow}

\newcommand{\til}{\rightarrow}

\newcommand {\del}{ \; \big| \;}

\newcommand {\mgdv}{\big\{}

\newcommand {\mgdh}{\big\}}

\newcommand {\Intv}{\Big[}
\newcommand {\intv}{\big[}

\newcommand {\Inth}{\Big]}
\newcommand {\inth}{\big]}

\newcommand {\for}{\bigcup}

\newcommand {\snit}{\bigcap}

\newcommand {\og}{\; \land \;}
\newcommand {\eller}{\; \vee\;}
\newcommand{\ikke}{\lnot}

\newcommand {\go} {\mathfrak}
\newcommand {\ku} {\mathcal}
\newcommand {\un} {\underline}

\newcommand {\e} {\exists}
\renewcommand {\a} {\forall}

\newcommand{\fed}{\boldsymbol}

\newcommand{\pf}{

\smallskip

\noindent {\it Proof : }}

\newcommand{\pff}{$\hfill  \Box$

\smallskip }

\newcommand{\CA}{\fed{\Pi}^1_1}
\newcommand{\PCA}{\fed{\Sigma}^1_2}

\newcommand{\bb}{bb_\Q(e_i)}

\newtheorem{thm}{Theorem}

\newtheorem{lemme}[thm]{Lemma}
\newtheorem{prop}[thm]{Proposition}
\newtheorem{defi}[thm]{Definition}

\newtheorem{ex}[thm]{Example}

\begin{document}
\maketitle
\begin{abstract}
A Banach space  contains either a minimal subspace or a continuum
of incomparable subspaces. General structure results for analytic
equivalence relations are applied in the context of Banach spaces
to show that if $E_0$ does not reduce to isomorphism of the
subspaces of a space, in particular, if the subspaces of the space
admit a classification up to isomorphism by real numbers, then any
subspace with an unconditional basis is isomorphic to its square
and hyperplanes and has an isomorphically homogeneous subsequence.
\end{abstract}

\section{Introduction.}

This paper contains results in the intersection of the geometry of
Banach spaces and descriptive set theory. The general problem of
our  study is a generalisation of the homogeneous space problem.
Namely, what can be said about a Banach space with ``few'' non
isomorphic subspaces? In particular, will such a space necessarily
satisfy more regularity properties than a general space? Will it
necessarily have subspaces of a given type?

The paper is divided into two parts, of which the first contains a
proof of the following:

\begin{thm}\label{incomparable} Let $X$ be an infinite dimensional
Banach space. Then $X$ contains either a minimal subspace or a
continuum of pairwise incomparable subspaces.\end{thm}

Recall that two spaces are said to be incomparable if neither of
them embed into the other, and a space is minimal if it embeds
into all of its infinite dimensional subspaces.

The homogeneous space problem, which was solved in the positive by
the combined efforts of Gowers \cite{g}, Komorowski and
Tomczak-Jaegermann \cite{kj}, is the problem of whether any
infinite dimensional space, isomorphic to all its infinite
dimensional subspaces, must necessarily be isomorphic to $\ell_2$.
As a continuation of this one can ask how many isomorphism classes
of subspaces  a non Hilbertian space has to contain. Infinitely
many? A continuum? Even for some of the classical spaces this
question is still open, though recent progress has been made by
Ferenczi and Galego \cite{fg}.

Our theorem and proof turn out to have something to say about the
following two problems of Gowers. (\cite{g}, Problems 7.9 and
7.10):
\begin{itemize}
\item Determine which  partial orders that can be realised as
the set of subspaces of an infinite dimensional Banach space under
the relation of embeddability. Or at least find strong conditions
such a partial order must necessarily satisfy.
\item Find further applications of the main determinacy result in
\cite{g}. In particular, are there any applications that need its
full strength, i.e., that need it to hold for analytic and not
just open sets?\end{itemize}

Our Theorem \ref{incomparable} says that any such partial order
must either have a minimal element or an antichain of continuum
size. And, as will be evident, the proof does in fact very much
need the full strength of the determinacy result.

We mention that our proof relies heavily on methods of logic and
we have therefore included a short review of the most basic
notions of set theory indispensable to understand the proof. Also
for the benefit of the non analyst we recall some standard notions
from Banach space theory.

Before presenting the results of the second part we will first
need this brief review.

\subsection{Descriptive set theory.}

Our general reference for descriptive set theory will be the book
by Kechris, \cite{ke}, whose notation will be adopted here.

A {\em Polish} space is a separable completely metrisable space. A
measurable space, whose algebra of measurable sets are the Borel
sets of some Polish topology, is said to be {\em standard Borel}.
These spaces turn out to be completely classified up to Borel
isomorphism by their cardinality, that can either be countable or
equal to that of the continuum. A subset of a standard Borel space
is {\em analytic} if it is the image by a Borel function of some
standard Borel space and {\em coanalytic} if its complement is so.
It is {\em $C$-measurable} if it belongs to the smallest
$\sigma$-algebra containing the Borel sets and closed under the
Souslin operation. In particular, analytic sets are $C$-measurable
as they can be obtained by the Souslin operation applied to a
sequence of Borel sets. $C$-measurable sets in Polish spaces
satisfy most of the classical regularity properties, such as
universal measurability and the Baire property. We denote by
${\fed{\Sigma}^1_1}$, $\CA$ and $\PCA$ the classes of analytic,
coanalytic and Borel images of coanalytic sets respectively. A
classic result of Sierpinski states that any $\PCA$ set is the
union of $\aleph_1$ Borel sets.

Let $X$ be a Polish space and $\ku F(X)$ denote the set of closed
subsets of $X$. We endow $\ku F(X)$ with the following
$\sigma$-algebra that renders it a standard Borel space. The
generators are the following sets, where $U$ varies over the open
subsets of $X$:
$$
\mgdv F \in \ku F (X) \del F \cap U \neq \tom \mgdh
$$
The resulting measurable space is called the Effros Borel space of
$X$.

Fix some basis $\{U_n\}$ for the space $C(\ca)$ and define the
Borel set $\go B$ by:
$$
\go B=\{F\in \ku F(C(\ca))\del \a n\; (0\in U_n\til F\cap U_n \neq
\tom)\og\a n,m,l\; \a r,t\in \Q
$$
$$
(rU_n+tU_m\subseteq U_l\og F\cap U_n\neq \tom\og F\cap U_m\neq\tom\til
F\cap U_l\neq \tom)\}
$$
This evidently consists of all the closed linear subspaces of
$C(\ca)$ and, as $C(\ca)$ is isometrically universal for separable
Banach spaces, any separable Banach space has an isometric copy in
$\go B$. We can therefore view $\go B$ as the standard Borel space
of all separable Banach spaces. When one wants to restrict the
attention to the subspaces of some particular space $X$ one only
needs to consider the Borel subset $\{Y\in \go B\del Y\subseteq
X\}$. Moreover, it is not hard to see that most reasonably
definable properties and relations are $\PCA$ in $\go B$ or $\go
B^n$; for example, the relations of isometry and isomorphism are
both analytic in $\go B^2$ exactly as expected.

A theme of descriptive set theory, that has been extensively
developed the last fifteen years or so, is the Borel reducibility
ordering of analytic equivalence relations on standard Borel
spaces.

This ordering is defined as follows: Suppose $E\subset X^2$ and
$F\subset Y^2$ are analytic equivalence relations on standard
Borel spaces $X$ and $Y$. We say that $E$ is {\em Borel reducible}
to $F$, in symbols $E\leq_BF$, if there is a Borel measurable
function $f:X\til Y$ such that for all $x,y\in X$:
$$xEy\equi f(x)Ff(y)$$
Moreover, when $X$ and $Y$ are Polish and $f$ can be taken to be
continuous, we write $E\leq_cF$.

Heuristically, $X$ represents a class of mathematical objects
(e.g., separable Banach spaces) that we wish to classify up to
$E$-equivalence (e.g., isomorphism) by complete invariants
belonging to some other category of mathematical objects. A
reduction $f\colon X\til Y$ of $E$ to $F$ corresponds then to a
classification of $X$-objects up to $E$-equivalence by $Y$-objects
up to $F$-equivalence.

Another way of viewing the Borel reducibility ordering is as a
refinement of the concept of cardinality. It provides a concept of
relative cardinality for quotient spaces in the absence of the
axiom of choice. For a reduction of $E$ to $F$ is essentially an
injection of $X/E$ into $Y/F$ admitting a Borel lifting from $X$
to $Y$.

A few words on the power of the continuum: We say that an analytic
equivalence relation $E$ on a standard Borel space $X$ has a {\em
continuum of classes} if there is an uncountable Borel set
$B\subset X$ consisting of pairwise $E$-inequivalent points. This
is known to be stronger than just demanding that there should be
some bijection between the set of classes  and $\R$. The are for
example analytic equivalence relations having exactly $\aleph_1$
many classes, but not having a continuum of classes (in the above
sense) in any model of set theory. But an uncountable Borel set is
always Borel isomorphic to $\R$, independently of the size of the
continuum.

If $A$ is some infinite subset of $\N$, we denote by $[A]^\N$ the
space of all infinite subsets of $A$ equipped with the topology
induced by the product topology on $2^A$. Furthermore, for two
sets $A$ and $B$ we write $A\subset^*B$ iff $A\setminus B$ is
finite. Then $A\subsetneq^*B$ iff  $A\subset^*B$ but
$B\not\subset^*A$. Also, when $A\subset \N$ and $k\in \N$ we let
$A/k=\{n\in A\del n>k\}$. We will occasionally also consider
natural numbers as ordinals, so that $n=\{0,1,\ldots,n-1\}$.

We will repeatedly use the following result of Ellentuck extending
results of Galvin-Prikry for Borel sets and Silver for analytic
sets: if $\ku A\subset[\N]^\N$ is a C-measurable set, then there
is some $A\in [\N]^\N$ with either $[A]^\N \subset \ku A$ or
$[A]^\N \cap \ku A=\tom$.

This has the consequence that if $f\colon \ram\til X$ is some
$C$-measurable function with values in some Polish space $X$, then
there is some $A\in\ram$ such that $f$'s restriction to $[A]^\N$
is continuous.

Among the simpler analytic equivalence relations are those that
admit a classification by real numbers, i.e., those that are Borel
reducible to the identity relation on $\R$. These are said to be
{\em smooth}. It turns out that among \un {Borel} equivalence
relations there is a minimum, with respect to $\leq_B$, non smooth
one, which we denote by $E_0$ (see \cite{hkl}). It is defined on
$\ram$ as the relation of eventual agreement, i.e.:
$$AE_0B\equi \e n\; A/n=B/n$$
To see that $E_0$ is non smooth, suppose towards a contradiction
that $f\colon \ram \til \R$ is a Borel function such that
$AE_0B\equi f(A)=f(B)$. Then there is some infinite $C\subset \N$
such that the restriction of $f$ to $[C]^\N$ is continuous. But,
as the equivalence class of $C$ is dense in $[C]^\N$, this means
that $f$ is constant on $[C]^\N$, contradicting that $[C]^\N$
intersects more than one equivalence class.

On the other hand any uncountable Borel set $B\subset \ram $ of
pairwise almost disjoint sets will witness that $E_0$ has a
continuum of classes.

From this it follows that any analytic equivalence relation to
which $E_0$ reduces has a continuum of classes, but does not admit
a classification by real numbers.

After these preliminary remarks we can state our second result.

\begin{thm} Let $X$ be a Banach space with an unconditional basis
$(e_n)$. If $E_0$ does not Borel reduce to isomorphism between
subspaces generated by subsequences of the basis (and in
particular if these admit a classification by real numbers), then
any space spanned by a subsequence is isomorphic to its square and
hyperplanes. Furthermore, there is a subsequence of the basis such
that all of its subsequences span isomorphic spaces.\end{thm}

For example, as the usual basis of Tsirelson's space does not have
a subsequence all of whose subsequences span isomorphic spaces,
this shows that there is no isomorphic classification of the
subspaces of Tsirelson's space by real numbers.

This result can be coupled with Gowers' dichotomy \cite{g}
proving:

\begin{thm}  Let $X$ be a separable Banach space. Either $E_0$ Borel
reduces to isomorphism between its subspaces or $X$ contains a
reflexive subspace with an unconditional basis all of whose
subsequences span isomorphic spaces.\end{thm}

For the above we will need some Ramsey type results for product
spaces and some constructions for reducing $E_0$. These results
seem to have an independent interest apart from their applications
to Banach space theory in that they classify minimal counter
examples to Ramsey properties in product spaces. Let us just state
one of these:

\begin{thm}Let $E$ be an analytic equivalence relation on $\ram$
invariant under finite changes. Either $E_0$ Borel reduces to $E$
or $E$ admits a homogeneous set.\end{thm}

\subsection{Schauder bases.}
Let $X$ be some separable Banach space and $( e_i)$ a non zero
sequence in $X$. We say that  $( e_i )$ is a {\em basis} for $X$
if any vector $x$ in $X$ can be uniquely written as a norm
convergent series $x=\sum a_ie_i$. In that case, the biorthogonal
functionals $e^*_k(\sum a_ie_i):=a_k$ and the projections
$P_n(\sum a_ie_i):=\sum_{i=0}^n a_ie_i$ are in fact continuous and
moreover their norms are uniformly bounded.

If $( e_i )$ is some non zero sequence that is a basis for its
closed linear span, written $\intv e_i \inth$, we say that it is a
{\em basic sequence} in $X$. The property of $(e_i)$ being a basic
sequence can also equivalently be stated as the existence of a
constant $K\geq 1$ such that for any $n\leq m$ and
$a_0,a_1,\ldots,a_m\in \R$:
$$\|\sum_{i=0}^n a_ie_i\|\leq K\|\sum_{i=0}^m a_ie_i\|$$

Suppose furthermore that for any $x=\sum a_ie_i$ the series
actually converges unconditionally, i.e., for any permutation
$\sigma$ of $\N$ the series $\sum a_{\sigma(i)}e_{\sigma(i)}$
converges to $x$. Then the basic sequence is said to be {\em
unconditional}.

Again, being an unconditional basis for some closed subspace
(which will be denoted by `unconditional basic sequence') is
equivalent to there being a constant $K\geq 1$, such that for all
n, $A\subset \mgdv 0,\ldots,n\mgdh$ and $a_0,\ldots,a_n\in \R$
$$\|\sum_{i\in A} a_ie_i\|\leq K\|\sum_{i=0}^n a_ie_i\|$$
We will in general only work with {\em normalised} basic
sequences, i.e., $\|e_i\|\equiv 1$, which always can be obtained
by taking $e'_i:=\frac{e_i}{\|e_i\|}$.

Given some vector $x\in span(e_i)$ let its {\em support},
$supp(x)$, be the set of indices $i$ with $e^*_i(x)\neq 0$. For
$k\in \N$ and $x,y\in span(e_i)$ we write $k<x$ if $k<\min
supp(x)$ and $x<y$ if $\max supp(x)<\min supp(y)$. A {\em block
basis}, $(x_i)$, over a basis $(e_i)$ is a finite or infinite
sequence of vectors in $span(e_i)$ with $x_0<x_1<x_2<\ldots$. This
sequence will also be basic and in fact unconditional in case
$(e_i)$ is so.

Two basic sequences $(e_i)$ and $(t_i)$ are called {\em
equivalent}, in symbols $(e_i)\approx (t_i)$,  provided a series
$\sum a_ie_i$ converges if and only if $\sum a_it_i$ converges.
This can also be stated as saying that $T\colon e_i\mapsto t_i$
extends to an invertible linear operator between $[e_i]$ and
$[t_i]$. The quantity $\|T\|\cdot\|T^{-1}\|$ is then the constant
of equivalence between the two bases.

A basis that is equivalent to all of its subsequences is said to
be {\em subsymmetric}. A simple diagonalisation argument then
shows that it must be uniformly equivalent to all of its
subsequences.

Two basic sequences $(e_i)$ and $(t_i)$ are said to be {\em
permutatively equivalent} if there is some permutation $\sigma$ of
$\N$ such that $(e_i)$ and $(t_{\sigma(i)})$ are equivalent.

\section{Incomparable and minimal subspaces.}

Two Banach spaces $X$ and $Y$ are called incomparable in case
neither of them embed isomorphically into the other. $X$ is said
to be minimal if it embeds into all of its infinitely dimensional
subspaces and $X$ itself is infinite dimensional.

Our proof of the first theorem will proceed by a reduction to an
analysis of Borel partial orders due to L. Harrington, D. Marker
and S. Shelah (see \cite {hms}). Instrumental in our reduction
will be the determinacy result of Gowers on certain games in
Banach spaces (see \cite{g}), which will guarantee that some
choices can be done uniformly, a fact that is needed for
definability purposes. Moreover, we will use some ideas of J.
Lopez-Abad on coding reals with inevitable subsets of the unit
sphere of a Banach space (see \cite {la}).

We mention that it was shown by a simpler argument in \cite {fr2}
by V. Ferenczi and the author that any Banach spaces either
contains a minimal subspace or a continuum of non isomorphic
subspaces.

For facility of notation we write $X\emb Y$ if $X$ embeds
isomorphically into $Y$ and will always suppose the spaces we are
working with to be separable infinite dimensional. Then $\emb$
restricted to the standard Borel space of subspaces of some
separable Banach space becomes an analytic quasi-order, i.e.,
transitive and reflexive. So the result above amounts to saying
that either $\emb$ has a minimal element or a perfect antichain.

Suppose $(e_i)$ is a normalised basic sequence with norm denoted
by $\|\cdot\|$. We call a normalised block vector $x$ with finite
support rational if it is a scalar multiple of a finite linear
combination of $(e_i)$ with rational coordinates. Notice that
there are only countably many rational (finite) block vectors,
which we can gather in a set $Q$ and give it the discrete
topology. Let $bb_\Q(e_i)$ be the set of block bases of $(e_i)$
consisting of rational normalised block vectors, which is easily
seen to be a closed subspace of $Q^\N$, which is itself a Polish
space. Moreover the canonical function sending $X\in\bb$ to its
closed span in $\go B$ is Borel, so the relations of isomorphism,
etc., become analytic on $\bb$.

We recall the following classical facts: Any  infinite dimensional
Banach space contains an infinite normalised basic sequence
$(e_i)$. Moreover, if $Y$ is any subspace of $[e_i]$, then it
contains an isomorphic perturbation of a block basic sequence of
$(e_i)$. Again any block basic sequence is equivalent to some
member of $\bb$. So this explains why we can concentrate on $\bb$
if we are only looking for minimal subspaces.

For $X,Y\in bb_\Q(e_i)$, let $X\leq Y$ if $X$ is a blocking of
$Y$, i.e., if any element of $X$ is a linear combination over $Y$.
Note that this does not imply that they are rational block vectors
over $Y$, but only over $(e_i)$. Moreover, if $Y=(y_i),X=(x_i)\in
bb_\Q(e_i)$, put $Y\leq^*X$ if for some $k$, $(y_i)_{i\geq k}\leq
X$. Also, for $\Delta=(\delta_i)$ an infinite sequence of strictly
positive reals write $d(X,Y)<\Delta$ if $\a i\;
\|x_i-y_i\|<\delta_i$.

Put $X\approx Y$ if the bases are equivalent and $X\iso Y$ if they
span isomorphic spaces. Then a classical perturbation argument
shows that there is some $\Delta$ depending only on the constant
of the basis, such that for any $X,Y\in bb_\Q(e_i)$ if
$d(X,Y)<\Delta$, then $X\approx Y$ and in particular $X\iso Y$.
Put also $X=(x_i)\backsimeq Y=(y_i)$ if $\e k\; \a i\geq k\;
x_i=y_i$. Then evidently $X\backsimeq Y$ implies $X\approx Y$.

For a subset $A\subset bb_\Q(e_i)$ let
$A^*=\{Y\in bb_\Q(e_i)\del \e X\in A \; X\backsimeq Y\}$ and
$A_\Delta=\{Y\in bb_\Q(e_i)\del \e X\in A \; d(X,Y)<\Delta\}$.
Notice that if $A$ is analytic so are both $A^*$ and $A_\Delta$.
Again $[Y]=\{X\in bb_\Q(e_i)\del X\leq Y\}$.
Such an $A$ is said to be large in $[Y]$ if for any $X\in [Y]$ we have
$[X]\cap A\neq \tom$.

For $(e_i)$ a given normalised basis, $A\subset bb_\Q(e_i)$ and
$X\in bb_\Q(e_i)$, the Gowers game  $\Game^A_X$ is defined as
follows: Player I plays in the $k$'th move of the game a rational
normalised block vector $y_k$ of $(e_i)$ such that $y_{k-1}<y_k$
and $y_k$ is a block on $X$. Player II responds by either doing
nothing or playing a rational normalised block vector $x$ such
that $x\in [y_{l+1}, \ldots, y_k]$ where $l$ was the last move
where II played a vector. So player II wins the game if in the end
she has produced an infinite rational block basis $X=(x_i)\in A$.
This is an equivalent formulation due to J. Bagaria and J.
Lopez-Abad (see \cite {bla}) of Gowers' original game.

Gowers \cite{g} proved that if $A\subset \bb $ is analytic, large
in $[Y]$ and $\Delta$ is given, then for some $X\in [Y]$ II has a
winning strategy in the game $\Game^{A_\Delta}_X$.

We mention also a result of Odell and Schlumprecht \cite{os}
obtained from their solution to the distortion problem: If $E$ is
an infinite dimensional Banach space not containing $c_0$, there
are an infinite dimensional subspace $F$ and $A,B\subset \ku S_F$
of positive distance such that any infinite dimensional subspace
of $F$ intersects both $A$ and $B$.

The following was shown in \cite{fr2}:

\begin {lemme}$(\fed{MA})$ Let $A\subset \bb$ be linearly ordered under
$\leq^* $ of cardinality strictly less than the continuum. Then
there is some $X\in\bb$ such that $X\leq^*Y$ for all $Y\in A$.\end
{lemme}

From this lemma one gets the following:

\begin{lemme}\label {PCA determinacy}$(\fed{MA+\neg CH})$ Suppose
$W\subset \bb$ is a $\fed{\Sigma}^1_2$ set, large in some $[Y]$
and  $\Delta>0$. Then II has a winning strategy in
$\Game^{W^*_\Delta}_X$ for some $X\in[Y]$.\end{lemme}

\pf Let $W=\for_{\omega_1}V_\xi$ be a decomposition of $W$ as an
increasing union of $\aleph_1$ Borel sets. We claim that some
$V^*_\xi$ is large in $[Z]$ for some $Z\in [Y]$, which by Gowers'
theorem will be enough to prove the lemma. So suppose not and find
$Y_0\in [Y]$ such that $[Y_0]\cap V_0^*=\tom$. Repeating the same
process and diagonalising at limits, we find $Y_\xi\in [Y]$ for
$\xi<\omega_1$ such that $[Y_\xi]\cap V_\xi^*=\tom$ and
$Y_\xi\leq^*Y_\eta$ for $\eta<\xi$. By the above lemma there is
some $Y_\infty=(y_i)\in [Y]$ with $Y_\infty\leq^*Y_\xi$ for all
$\xi<\omega_1$.

 We claim that $[(y_{2i})]\cap W=\tom$. Otherwise, for
 $Z=(z_i)\in [(y_{2i})]\cap W$
 find $\xi<\omega_1$ such that $Z\in V_\xi$. Now as $(y_i)\leq ^*Y_\xi$
 there is some $k$
 with $(y_i)_{i\geq k}\leq Y_\xi$, but then
$$(y_{2i})\backsimeq(y_k,y_{k+1},y_{k+2},\ldots,y_{2k-1},y_{2k},
y_{2(k+1)},y_{2(k+2)},
\ldots)\leq Y_\xi$$
 One now easily sees that there is some $(x_i)$ with
$$(z_i)\backsimeq(x_i)\leq  (y_k,y_{k+1},y_{k+2},\ldots,y_{2k-1},
y_{2k},y_{2(k+1)},
y_{2(k+2)},\ldots)$$
 whereby $(x_i)\in V_\xi^*$ contradicting $V_\xi^*\cap [Y_\xi]=\tom$.

Therefore $[(y_{2i})]\cap W=\tom$, again contradicting the
largeness of $W$.\pff

\begin {lemme}\label{non minimality}$(\fed{MA+\neg CH})$ Suppose that
$(e_i)$ is a basic sequence such that $[e_i]$ does not contain a
minimal subspace. Then for any $Y\in\bb$ there are a $Z\in [Y]$
and a Borel function $g\colon [Z]\til [Z]$, with $g(X)\leq X$ and
$X\not\sqsubseteq g(X)$ for all $X\in [Z]$.\end{lemme}

\pf As $c_0$ is minimal, $[e_i]$ does not contain $c_0$.
Therefore, by the solution to the distortion problem by Odell and
Schlumprecht, we can by replacing $(e_i)$ by a block suppose that
we have  two positively separated sets $F_0, F_1$ of the unit
sphere, such that for any $X\in \bb$ there are rational normalised
blocks $x,y$ on $X$ with $x\in F_0$ and $y\in F_1$. We call such
sets inevitable.

Let $D=\{X=(x_i)\in \bb\del \a i\; x_i\in F_0\cup F_1\}$ and for
$X\in D$ let $\alpha(X)\in\ca$ be defined by $\alpha(X)(i)=0\hviss
x_i\in F_0$. Then $D$ is easily seen to be a closed subset of
$\bb$ and $\alpha\colon D\til \ca$ to be continuous. Furthermore
by the inevitability of $F_0$ and $F_1$ we have that $D$ is large
in every $[Y]$.

Let $\Q_*^{<\N}$ be the set of finite non identically zero
sequences of rational numbers given the discrete topology. Then
$(\Q_*^{<\N})^\N$ is Polish. Define for any $Y\in\bb$ and
$(\un{\lambda}_i)\in(\Q_*^{<\N})^\N$ the block basis
$(\un{\lambda}_i)\cdot Y$ of $Y$ in the obvious way, by taking the
linear combinations given by $(\un{\lambda}_i)$.

Fix also some perfect set $P$ of almost disjoint subsets of $\N$
seen as a subset of $\ca$ and let $\beta\colon
P\hviss(\Q_*^{<\N})^\N$ be a Borel isomorphism.

Again $E=\{X\in D\del \alpha(X)\in P\}$ in large and closed in $\bb$.

Then the set
$$W=\{X=(x_i)\in  \bb\del (x_{2i})\in E\og (x_{2i+1})\not
\sqsubseteq\beta\circ\alpha
((x_{2i}))\cdot(x_{2i+1})\}$$
 is coanalytic. We claim moreover that it is large in $\bb$.

To see this, let $Y\in \bb$ be given and take by inevitability of
$F_0$ and $F_1$ some $(z_i)\in [Y]$ with $z_{3i}\in F_0$ and
$z_{3i+1}\in F_1$. As  $[z_{3i+2}]$ is not minimal there is some
$X\leq(z_{3i+2})$ such that $(z_{3i+2})\not\sqsubseteq X$. Take
some $(\un{\lambda}_i)\in(\Q_*^{<\N})^\N$ such that
$(\un{\lambda}_i)\cdot(z_{3i+2})\approx X$ and
$(z_{3i+2})\not\sqsubseteq (\un{\lambda}_i)\cdot(z_{3i+2})$. We
can now define some $(v_i)$ such that either $v_{2i}=z_{3i}$ or
$v_{2i}=z_{3i+1}$, $\beta\circ\alpha((v_{2i}))=(\un{\lambda}_i)$
and  $v_{2i+1}=z_{3i+2}$. This ensures that $(v_i)\in W$. So as
$(v_i)\leq(z_i)$ it is in $[Y]$ and $W$ is indeed large.

Take now some $\Delta=(\delta_i)$ depending on the basic constant
as above with $\delta_i<\frac12d(F_0,F_1)$. By the preceding lemma
we can find a $Y\in \bb$ such that II has a winning strategy
$\sigma$ in the game $\Game^{W^*_\Delta}_Y$.

Suppose that $=(x_i)$ has been played by II according to the
strategy $\sigma$ as a response to $Z$ played by I. As $\sigma$ is
winning, $X\in W^*_\Delta$. Define $\gamma(X)\in\ca$ by
$\gamma(X)(i)=0$ if $d(x_{2i},F_0)<\delta_{2i}$ and
$\gamma(X)(i)=1$ otherwise. Then $\gamma$ is Borel from
$W^*_\Delta$ to $\ca$, and furthermore there is a unique
$\gamma^*(X)\in P$ such that $\e k\; \a i\geq k\;
\gamma(X)(i)=\gamma^*(X)(i)$. This is because $P$ was chosen to
consist of almost disjoint subsets of $\N$. Again $X\mapsto
\gamma^*(X)$ is Borel.

Take some $U=(u_i)\in W$ such that $\a^\infty n\;
\|u_n-x_n\|<\delta_n$. Then $\alpha(U)=\gamma^*(X)$,
$(u_{2i+1})\approx (x_{2i+1})$ and $(u_{2i+1})\not\sqsubseteq
\beta\circ\alpha(U)\cdot(u_{2i+1})$. So due to the equivalence
invariance of the basis by $\Delta$ perturbations we have
$(x_{2i+1}) \not\sqsubseteq\beta\circ\gamma^*(X)\cdot(x_{2i+1})$.

Let $V\in [X]$ be the normalisation of
$\beta\circ\gamma^*(X)\cdot(x_{2i+1})$. The function $g\colon
Z\mapsto V$ is Borel and obviously

$$V\approx\beta\circ\gamma^*(X)\cdot(x_{2i+1})\leq (x_{2i+1})\leq
Z$$

and as $$(x_{2i+1})
\not\sqsubseteq\beta\circ\gamma^*(X)\cdot(x_{2i+1})$$

also $Z\not\sqsubseteq V$.

\pff

A Banach space is called quasi-minimal if any two subspaces have
further isomorphic subspaces. The following is a standard
observation.

\begin{lemme}Suppose $[e_i]$ is quasi-minimal. Then $\sqsubseteq$ is
downwards $\sigma$-directed on $\bb$, i.e., any countable family
has a common minorant.\end{lemme}

\pf Suppose that $Y_i\in\bb$ are given, then define inductively
$Z_i\in[Y_0]$ such that $Z_i\sqsubseteq Y_i$ and $Z_{i+1}\leq
Z_i$. Take some $Z=(z_i)\leq^* Z_n$ for all $n$ and notice as in
the proof of lemma \ref {PCA determinacy} that
$(z_{2i})\sqsubseteq Z_n$ for all $n$.\pff

\begin{lemme}If $R$ is a downwards $\sigma$-directed Borel quasi-order
on a standard Borel space $X$. Then either $R$ has a perfect
antichain or a minimal element.\end{lemme}

\pf This is a simple consequence of the results of L. Harrington,
D. Marker and S. Shelah \cite{hms}, as we will see. Suppose that
$R$ did not have a perfect antichain, then by their results there
is a countable partition $X=\for X_n$ into Borel sets, so that $R$
is total on each piece, i.e., $R$ can be written as a countable
union of $R$-chains.

Applying another of their results this implies that for some
countable ordinal $\alpha$ there are Borel functions $f_n\colon
X_n\til 2^{\alpha}$, such that for any $x,y\in X_n$:
$$yRx\equi x\leq_{lex}y$$
Where $\leq_{lex}$ is the usual lexicographical ordering. In their
terminology, $R$ is linearisable on each $X_n$.

One can easily check that any subset of $2^{\alpha}$ has a
countable subset cofinal with respect to $\leq_{lex}$, so pulling
it back by $f_n$ it becomes coinitial in $R\begr_{X_n}$. Putting
all these sets together one gets a countable subset of $X$
coinitial with respect to $R$. So by downwards
$\sigma$-directedness there is therefore a minimal element in
$X$.\pff

After this series of lemmas we can now prove the theorem:

\begin{thm}Let $X$ be an infinite dimensional Banach space. Then $X$
contains either a minimal subspace or a continuum of pairwise
incomparable subspaces.\end{thm}

\pf By Gowers' quadrichotomy $X$  contains either a quasi-minimal
subspace or a subspace with a basis such that any two disjointly
supported subspaces are totally incomparable (see Gowers \cite{g}
theorem 7.2 and the fact that H.I. spaces are quasi-minimal). In
the latter case any perfect set of almost disjoint subsets of $\N$
will give rise to subsequences of the basis spanning totally
incomparable spaces, which would prove the theorem. So we can
suppose that $X=[e_i]$ is quasi-minimal for some basis $(e_i)$. If
$X$ does not contain a minimal subspace, we can choose $Z\in\bb$
and the Borel function as above (under $\fed{MA+\neg CH}$ of
course). So define the following property on subsets $A,B$ of
$[Z]^2$:
$$\Phi(A,B)\Leftrightarrow $$
$$\a Y,V,W\in[Z]\;\intv (Y,V)\notin A\eller (V,W)\notin A\eller
(Y,W)\notin B
\inth\og\a Y\in [Z]\; (Y,g(Y))\notin A$$

We see that $\Phi$ is $\fed\Pi^1_1$ on $\fed\Sigma^1_1$,
hereditary and continuous upwards in the second variable.
Furthermore, $\Phi(\sqsubseteq,\not\sqsubseteq)$, so by the second
reflection theorem (see Kechris \cite{k} theorem (35.16)) there is
some Borel set $R$ containing $\sqsubseteq$ such that
$\Phi(R,\compl R)$. But then $R$ is a Borel quasi-order, downwards
$\sigma$-directed, as it contains $\sqsubseteq$, and without a
minimal element, as witnessed by $g$. So $R$ has a perfect
anti-chain by the previous lemma, which then is an antichain for
$\sqsubseteq$ too.

The statement is therefore proved under the additional hypothesis
of Martin's axiom and the negation of the continuum hypothesis. We
will see that this is in fact sufficient to prove the theorem. By
standard metamathematical facts and Shoenfield's absoluteness
theorem it is enough to show that the statement we wish to prove
is $\Sigma^1_2$.

It was proved by Ferenczi and the author in \cite{fr2} that the
property of having a block minimal subspace was $\Sigma^1_2$. For
using Gowers' determinacy result and codings as above, one can
continuously find an isomorphism between the space and a certain
subspace to testify the minimality. This proof can trivially be
modified to show that the property of having a minimal (i.e., not
necessarily block minimal) subspace is also $\Sigma^1_2$. For now
we only have to choose not a code for a subspace and an
isomorphism, but a code for a subspace and an embedding. For the
convenience of the reader, we have included the proof of this in
an appendix.

On the other hand, the property of having a perfect antichain is
obviously $\Sigma^1_2$ by just counting quantifiers. So these
remarks finish the proof.\pff

\section{Ramsey type results.}

We will show two Ramsey type results and afterwards some
applications to Banach space theory.

It is  well known that there are no nice Ramsey properties for the
product space $\ram\times\ram$ in contradistinction to the simple
Ramsey space $\ram$. That is, there are even quite simple
relations not admitting a square $[A]^\N\times[B]^\N$ that is
either included in or disjoint from the relation. An example of
this is the oscillation relation $\ku O$ defined by
$$
(a_n)\ku O(b_n)\equi
$$
$$
\e N\; \a n\;\Intv \#(k\;|\; a_n<b_k<a_{n+1}) \leq N \og
\#(k\;|\; b_n<a_k<b_{n+1})\leq N\Inth
$$
Where $\ram$ is seen as
the space of strictly increasing sequences of integers $(a_i)$.

The situation is very different if one replaces one of the factor
spaces by other Ramsey spaces and there are now very deep positive
theorems on so called polarised partition relations.

We are interested in the case when the relation on the product is
in fact a definable equivalence relation. Here the right question
seems to be when there is a cube $[A]^\N$ contained in one class.
Now if one lets two subsets of $\N$ be equivalent iff they have
the same minimal element, then the relation has exactly $\aleph_0$
classes and does not admit a homogeneous set.

On the other hand if the relation is invariant under finite
changes, such as $E_0$, then there are bigger chances that it
should have a homogeneous set. We will show that in the case of
analytic equivalence relations, $E_0$ is in fact the minimal
counterexample to the Ramsey property, in the sense that, if an
analytic equivalence relation is invariant under finite changes
and does not admit a homogeneous set, then it Borel reduces $E_0$.
In the same vein it is shown that if an analytic equivalence
relation does not admit a cube on which it has only countably many
classes, then it has at least a perfect set of classes. We notice
that both of these results are relatively direct consequences of
the Silver and Glimm-Effros dichotomies in the case of the
equivalence relation being Borel. But our results are motivated by
applications to isomorphism of separable Banach spaces, which is
true analytic, and the dichotomies are known not to hold in this
generality.

The following result was also found independently by S.
Todorcevic, albeit with a somewhat different proof:

\begin {thm}
Let $E$ be an analytic equivalence relation on $[\N]^\N$. Then
either $E$ has a continuum of classes or there is some $A\in
[\N]^\N$ such that $E$ only has a countable number of classes on
$[A]^\N$.

Moreover, $A$'s $E_0$ class will be a complete section for $E$ on
$[A]^\N$.
\end {thm}

\pf We will prove the theorem under $\fed{MA+\neg CH}$. By
Burgess' theorem (Exercise (35.21) in \cite{ke}) we can suppose
that $E$ has at most $\aleph_1$ classes $(C_\xi)_{\omega_1}$.
Define $P_\xi(A)\hviss [A]_{E_0}\cap C_\xi\neq \tom$ and notice
that this an analytic $E_0$-invariant property. We can by simple
diagonalisation find $(A_\xi)_{\omega_1}$, $A_\xi\subset^*A_\eta$
for $\eta<\xi<\omega_1$ such that $\a \xi<\omega_1$ either
$[A_\xi]^\N\subset P_\xi$ or $[A_\xi]^\N\subset \compl P_\xi$. And
by $\fed{MA+\neg CH}$ there is an $A\subset^*A_\xi,\a
\xi<\omega_1$.

Notice now that by $E_0$-invariance of $P_\xi$, if $B\subset^*A$
and $[A]^\N\subset P_\xi$ or   $[A]^\N\subset \compl P_\xi$ then
also $[B]^\N\subset P_\xi$, respectively   $[B]^\N\subset \compl
P_\xi$. So therefore $\a \xi <\omega_1\; [A]^\N\subset P_\xi$ or
$[A]^\N\subset \compl P_\xi$.

Suppose now that $B\in [A]^\N, B\in C_\xi$ then $P_\xi(B)$ and
therefore $[A]^\N\subset P_\xi$ and $P_\xi(A)$, i.e., $\e
A'E_0A\;\;\; A'EB$. This means that $[A]_{E_0}$ is a complete
section for $E$ on $[A]^\N$.

Let us now see that the statement of the theorem is absolute.
Saying that $E$ has a continuum of classes is equivalent to saying
that there is a compact perfect set $K\subset \ram$ consisting of
pairwise $E$-inequivalent points:
$$
\e K\subset \ram \;\textrm{compact, perfect}
\;\a x,y\in K (x=y \eller xEy)
$$
This is obviously a $\Sigma^1_2$ statement.

For the other case, notice that as $[A]_{E_0}$ is a complete
section for $E$ on  $[A]^\N$ there is by the Jankov-von Neumann
selection theorem a $C$-measurable selector  $f\colon [A]^\N\til
\N^\N$ and a Borel set $\ku B\subset \ram \times \ram\times \N^\N$
with $E=\pi_{\ram\times\ram}\ku B$, such that for $D\in [A]^\N$
there is an $A'E_0A$ with $\ku B\big ( D,A', f(D)\big)$. That is,
we can choose a witness to $D$ being $E$ equivalent to some
$A'E_0A$ in a $C$-measurable way. But any C-measurable function
can, using Ellentuck's theorem, be rendered continuous on a cube,
i.e., there is some $B\in [A]^\N$ such that $f$'s restriction to
$[B]^\N$ is continuous. So by the proof above the $E$-classes on
$[A]^\N$ are the same as the $E$-classes on $[B]^\N$ and the other
possibility can be written as:
$$\e B \;
\e f\colon [B]^\N\til \N^\N \;\textrm{continuous} \;\a D\in
[B]^\N\; \e B'E_0 B\; \;\ku B\big(D,B',f(D)\big)$$

This statement is $\Sigma_2^1$ as the quantifier $\e B'E_0 B$ is
over a countable set, so by Shoenfield absoluteness and standard
metamathematical facts it is enough to prove the result under
$\fed{MA+\neg CH}$.\pff

Our next results render explicit the connection with the Borel
reducibility ordering.

\begin{defi}
For $A, B\subset \N$ set $AE'_0B$ iff  $\e n \; |A\cap n| =|B\cap
n| \og A\setminus n= B\setminus n$
\end{defi}

It is easy to see that the equivalence class of any
infinite-coinfinite subset  of $\N$ is dense in $\ram$ and in fact
the equivalence relation is generically ergodic. Moreover, $E'_0$
is just a refinement of $E_0$.

\begin{lemme}
$E'_0$ is generically ergodic (i.e., any invariant set with the
Baire property  is either meagre or comeagre) and all classes
$[A]_{E'_0}$, for $A$ infinite-coinfinite, are dense.
\end{lemme}

\pf Since $\ram$ is cocountable in $\ca$ we can  restrict our
attention to it.  Suppose that some invariant set $\ku A$ is non
meagre, then there is some $a\subset [0,n]$ such that $\ku A$ is
comeagre in $D_{a,n}=\{A\in \ram\del A\cap [0,n]=a\}$. So for any
$D_{b,m}$ there are $c, d\in [0,k];\;  max(n,m)<k$ such that
$a\subset c,\; b\subset d,\; |c|=|d|$. Now for any $A\in
[\{k+1,k+2,\ldots\}]^\N$ we have $\phi(c\cup A):=(d\cup A)E'_0(
c\cup A)$ and $\phi$ is a homeomorphism of $D_{c,k}\subset
D_{a,n}$ with $D_{d,k}$. But that means that the image of $\ku A$
is comeagre in $D_{d,k}\subset D_{b,m}$ and is included in the
saturation of $\ku A$, which is $\ku A$. So $\ku A$ is comeagre in
the space.

If $A$ is infinite-coinfinite, then for any $D_{a,n}$ there are
$b,c\subset [0,k];\;b\supset a,\; n<k, \; A\cap [0,k]=c, \;  b\cap
[0,n]=a$ and $|b|=|c|$. So $A=(c\cup A/k)E'_0(b\cup A/k)\in
D_{b,k}\subset D_{a,n}$. And its class is dense. \pff

\begin{prop}
Let $E$ be a meagre equivalence relation on $\ca$ containing
$E'_0$. Then  $E_0\leq_B E$.
\end {prop}

\pf Let $(D_n)$ be a decreasing sequence of dense open sets, such
that $E\cap \snit_n D_n=\tom$.

We will inductively construct sequences $b_0^n, b_1^n\in \cantor$
for $n\in \N$ such that  for all $n$, $|b^n_0|=|b_1^n|,\;
\overline{\overline{b_0^n}}=\overline{\overline{b_1^n}}:=\#\{k\del
b_1^n(k)=1\}$. And if $a_s:={b^0_{s(0)}}^\smallfrown\ldots
^\smallfrown {b^{|s|-1}_{s(|s|-1)}}$ for all $s\in \cantor$, then
for any $s,t \in 2^n$, $N_{a_{s^\smallfrown 0}}\times
N_{a_{t^\smallfrown 1}}\subset D_{n+1}$.

Suppose that this can be done. Then define  $\alpha\mapsto \cup_n
a_{\alpha\begr n}=a_\alpha$. This is clearly continuous. If now
$\ikke \alpha E_0 \beta$ , then for infinitely many $n$,
$\alpha(n)\neq\beta(n)$. So for these $n\;$ $(a_\alpha,a_\beta)\in
N_{a_{\alpha\begr{n+ 1}}}\times N_{a_{\beta\begr n+1}}\subset
D_{n+1}$, which implies that $(a_{\alpha}, a_\beta)\in \snit_k D_k
\subset \compl E$.

Conversely, if $\alpha E_0 \beta$, then for some $N$, we have  $\a
n \geq N\; \alpha (n)=\beta (n)$. But then easily $a_\alpha=
{a_{\alpha\begr N}}^\smallfrown
 {b_{\alpha(N)}^N} ^\smallfrown {b_{\alpha (N+1)}^{N+1}}\ldots$ and
 $a_\beta={ a_{\beta\begr N}}^\smallfrown {b_{\alpha(N)}^N}^\smallfrown
 b_{\alpha (N+1)}^{N+1}\ldots$, so by the construction, $a_\alpha E'_0
 \beta$.

Now for the construction: Suppose that $b_0^n,b_1^n$ have been
chosen for $\a m<n$, enumerate  $2^n\times2^n$ by $(s_0,t_0),
\ldots, (s_k,t_k)$ and take $c_0^0,c_1^0\in \cantor$ such that
$N_{{a_{s_0}}^\smallfrown c_0^0}\times N_{{a_{t_0}}^\smallfrown
c_1^0}\subset D_n$. This can be done as $D_n$ is dense and open in
the product.

Prolong  $c_0^0,c_1^0$ to  $c_0^1,c_1^1$ respectively in such a way
that  $N_{{a_{s_1}}^\smallfrown c_0^1}\times N_{{a_{t_1}}^\smallfrown
c_1^1}\subset D_n$.

Again, prolong  $c_0^1,c_1^1$ to  $c_0^2,c_1^2$ respectively in
such a way that  $N_{{a_{s_2}}^\smallfrown c_0^2}\times
N_{{a_{t_2}}^\smallfrown c_1^2}\subset D_n$, etc.

Finally, prolong  $c_0^k,c_1^k$ to  $b_0^n,b_1^n$ respectively,
such that   $|b^n_0|=|b_1^n|,\;
\overline{\overline{b_0^n}}=\overline{\overline{b_1^n}}$. This
finishes the construction. \pff

For the following, we recall that $\a^* x \;R(x)$ means that the
set $\{x\del R(x)\}$ is comeagre, where $x$ varies over some
Polish space.

\begin{thm}\label{ramsey}
Let $E$ be an analytic equivalence relation on $\ram$ such that
$E'_0\subset E$, i.e.,  $E$  is $E'_0$-invariant. Then either
$E_0\leq_c E$ or there is some $A\in \ram $ such that $E$ only has
one class on $[A]^\N$.
\end {thm}

\pf By corollary 3.5 of \cite{HK}, if $E_0 \not\leq_c E$, then $E$
will be a decreasing  intersection of  $\aleph_1$ smooth
equivalence relations:

$$E=\snit_{\omega_1} E_\xi , \;\; E_\xi\subset E_\eta,\;\; \eta < \xi
< \omega_1$$

Let $f_\xi\colon \ram\til\R$ be a Borel reduction of $E_\xi $ to
identity on $\R$. Then for any $A\in \ram $, there is a $B\in
[A]^\N$ such that $f\begr [B]^\N$  is continuous. But since there
is a dense $E_\xi$-class the function has to be constant, that is,
there is only one class.

We construct inductively a $\subset^*$-decreasing sequence
$(A_\xi)_{\omega_1}$ of  infinite subsets of $\N$, with each
$A_\xi$ being homogeneous for $E_\xi$. Under $\fed{MA+\neg CH}$
such a sequence can be diagonalised to produce an infinite
$A_\infty\subsetneq^* A_\xi, \; \a \xi< \omega_1$. Now as
$A_\infty \subsetneq^* A_\xi$ it is easily seen that $A_\infty$ is
$E'_0$-equivalent with some subset of $A_\xi$ and therefore also
$E_\xi$-equivalent with $A_\xi$ itself. Furthermore, the same
holds for any infinite subset of $A_\infty$, so $A_\infty$ is
homogeneous for all of the $E_\xi$ and therefore for $E$ too.

As before one sees that the property of having a homogeneous set
is $\Sigma_2^1$,  so we need only check that continuously reducing
$E_0$ is $\Sigma^1_2$. But this can be written as:
$$\e f\colon \ram\til \ram \textrm{ continuous} \;\Big[ \a^*\alpha \in
\ram \;\a \beta E_0 \alpha  \;\a \gamma E_0 \alpha \;
f(\beta)Ef(\gamma)$$

$$\og \a \alpha , \beta \in \ram \big\{\alpha E_0 \beta \eller \ikke
f(\alpha)E f(\beta)\big\}\Big]$$

So as the quantifier $\a \beta E_0 \alpha$ is over a  countable
set and that the category quantifier $\a^*$ preserves analyticity
(see Theorem (29.22) in \cite{ke}), the statement is $\Sigma^1_2$.
\pff

\section{Applications to Banach space theory.}

Let $(e_i)$ be some basic sequence in a Banach space $X$ and
define  the following equivalence relation on $\ram$: $A\iso
B\equi[e_i]_A\iso[e_i]_B$. Then $\iso$ is analytic and extends
$E'_0$. For suppose that $AE'_0B$. Then $[e_i]_A$ and $[e_i]_B$
are spaces of the same finite codimension in $[e_i]_{A\for B}$ and
are therefore isomorphic. So, using the proposition, one sees that
if $ E_0 \not\leq_B\; \iso$, then $\iso$ must be non meagre and
therefore by Kuratowski-Ulam have a non meagre class, which again
by the lemma is comeagre.

To avoid trivialities, let us in the following suppose that all
Banach spaces considered are separable, infinite dimensional.

Gowers showed the following amazing result about the structure of
subspaces of a  Banach space: if $X$ is a Banach space, then it
contains either an unconditional basic sequence or an H.I.
subspace \cite{g}.

Here an H.I. (hereditarily indecomposable) space $Y$ is one in
which no two  infinite dimensional subspaces form a direct sum.
This property, which passes to subspaces, insures that $Y$ cannot
be isomorphic to any of its subspaces and cannot contain any
unconditional basic sequence. Therefore in the classification of
the subspaces of a Banach space one can always suppose to be
dealing with an H.I. space or a space with an unconditional basis.

\begin{prop}
Let $(e_i)$ be a basic sequence in a Banach space. Then either
$E_0$ Borel reduces  to isomorphism of spaces spanned by
subsequences of the basis or there will be some infinite $A\subset
\N$, such that for any infinite $B\subset A$: $[e_i]_A\iso
[e_i]_B$.
\end{prop}

\pf This follows from Theorem \ref{ramsey} as isomorphism is
$E'_0$-invariant.\pff

\begin{ex} Hereditarily Indecomposable spaces.
\end{ex}
Suppose that we are given a hereditarily indecomposable space $X$.
Then as any  Banach space contains a (conditional) basic sequence,
we can suppose that we have a basis $(e_i)$. By the above
proposition, if $E_0$ does not reduce, there would be a
subsequence spanning a space isomorphic to some proper subspace in
contradiction with the properties of H.I. spaces. So $E_0$ reduces
to isomorphism of its subspaces. The same reasoning shows, using
the first theorem, that it has a continuum of incomparable
subspaces.

\

A recent result due to Ferenczi and Galego \cite{fg} says that
$E_0$ Borel reduces to the isomorphism relation between subspaces
of $c_0$ and $\ell_1$. So if $E_0$ does not reduce to isomorphism
between the subspaces of an Banach space, then using Gowers'
dichotomy we can find a subspace with an unconditional basis.
Therefore by James' characterisation of reflexivity this basis
must span a reflexive space. All in all this gives us the
following:

\begin{thm}Let $X$ be an Banach space such that the isomorphism
relation between its subspaces does not reduce $E_0$. Then $X$
contains a reflexive subspace with an unconditional basis, all of
whose subsequences span isomorphic spaces.\end{thm}

Let us notice that if a basis $(e_i)$ has the property that no two
disjointly  supported block basic sequences are equivalent, then
one can easily show that this basis has the Casazza property and
moreover that it satisfies
$$
(e_i)_A\approx(e_i)_B\equi [e_i]_A\iso[e_i]_B\equi AE'_0B
$$
See the work of Gowers and Maurey, \cite{GM2},  for unconditional
examples of such bases. So as $E_0'$ and $E_0$ are Borel
bi-redicible, there are bases on which both equivalence and
isomorphism between subsequences are exactly of complexity $E_0$.

We will now see  an extension of some results by Ferenczi and the
author,  \cite{fr}, and Kalton, \cite{k}.

\begin{thm}
Let $(e_i)$ be an unconditional basic sequence. Then either $E_0$
Borel reduces  to isomorphism of spaces spanned by subsequences of
the basis or any space spanned by a subsequence is isomorphic to
its square and its hyperplanes. And there is some infinite
$A\subset \N$ such that for any infinite $B\subset A$,
$[e_i]_A\iso [e_i]_B$.
\end{thm}

\pf As before we can suppose we have some comeagre class $\ku
A\subset \ram$.   But then $\ku A$ is also comeagre in $\ca$ and
there is therefore a partition $A_0,A_1$ of $\N$ and subsets
$B_0\subset A_0, B_1\subset A_1$ such that for any $C\subset \N$,
if $C\cap A_0=B_0$ or $C\cap A_1=B_1$, then $C\in \ku A$. In
particular, $B_0,B_1,B_0\cup B_1\in \ku A$. Moreover, as the
complement operation is a homeomorphism of $\ca$ with itself,
there is some $C$ such that $C, \compl C\in \ku A$. So identifying
subsets of $\N$ with the Banach spaces they generate and using the
fact that the basis is unconditional, and therefore that disjoint
subsets form direct sums, we can calculate:
$$
\N=C\cup \compl C\iso C\oplus\compl C\iso B_0\oplus B_0\iso B_0
\oplus B_1\iso B_0 \cup B_1\iso B_0
$$
So $\N\in \ku A$ and $\ku A$ consists of spaces isomorphic to
their squares. Now for any $D\subset \N$:
$$
\N\oplus D\iso B_0\oplus B_1\oplus D\iso \intv B_0\cup(D\cap A_1)
\inth\oplus\intv B_1\cup (D\cap A_0)\inth\iso \N\oplus\N\iso\N
$$
This in particular shows that $[e_i]_\N$ is isomorphic to its
hyperplanes.

We notice now that the argument is quite general, in the sense
that we could have begun from any  $[e_i]_A$ instead of
$[e_i]_\N$, and therefore the results hold for any space spanned
by a subsequence. \pff

Kalton \cite{k} showed that in case an unconditional basis only has a
countable  number of isomorphism classes on the subsequences of the
basis, then the space spanned is isomorphic to its square and
hyperplanes.  The above result is along the same lines and we should
mention that one can get uniformity results with a bit of extra care in
the proof, see the article by Ferenczi and the author, \cite{fr}, for
this.

Notice that permutative equivalence between subsequences of a basis
induces an  analytic equivalence relation on $\ram$.

P. Casazza drew my attention to the following theorem from \cite{bclt}
(proposition 6.2).

\begin{thm}
{\emph{(Bourgain, Casazza, Lindenstrauss, Tzafriri)}} If
$(e_i)_\N$ is an  unconditional basic sequence permutatively
equivalent to all of its subsequences, then there is a permutation
$\pi$ of $\N$ such that $(e_{\pi(i)})_\N$ is subsymmetric.
\end{thm}
Their statement of the theorem is slightly more general, but the
general case  is easily seen to follow from the infinite dimensional
Ramsey theorem.
\begin{prop}
Let $(e_i)_\N$ is an unconditional basic sequence. Then either
$E_0$ reduces  to the relation of permutative equivalence of the
subsequences of the basis or there is some $A\in \ram$ such that
$(e_i)_A$ is subsymmetric.
\end{prop}

\pf Notice that $\sim_p$ (permutative equivalence) on $\ram$ is
$E'_0$-invariant, so  applying the Ramsey result we can suppose
that there is some $B\in\ram$ such that all $C\in[B]^\N$ are
$C\sim_p B$. Now there is some permutation $\pi$ of $B$ such that
$(e_{\pi(i)})_B$ is subsymmetric. Again choosing a strictly
increasing sequence $A=\{n_0,n_1,n_2,\ldots\}\subset B$ such that
$\pi(n_0)<\pi(n_1)<\pi(n_2)<\ldots$, we get a subsymmetric
$(e_i)_A$. \pff

{\em Acknowledgement} This article was written while the author
was a doctoral student in Paris under the direction of Alain
Louveau. I am sincerely grateful for all his help and the interest
he showed in my work.

\appendix\section*{Appendix}

For the convenience of the reader we include the proof to the
effect that having a minimal subspace is $\Sigma^1_2$. This is a
slightly amended version of the proof in \cite{fr2} showing that
having a block-minimal subspace is a $\Sigma^1_2$ property.

So suppose that $X$ has a minimal subspace and that it does not
contain $c_0$. Then it has a minimal subspace with a basis $(e_i)$
and positively separated inevitable sets $F_0$, $F_1$ in the unit
sphere of $[e_i]$. (Again this is by the results of Odell and
Schlumprecht.)

We let $D=\mgdv Y=(y_i)\in \bb\del \a i\; y_i\in F_0\cup F_1\mgdh$
and as in lemma \ref{non minimality} let $\alpha(Y)\in \ca$ be
defined by $\alpha(Y)(i)=0\hviss y_i\in F_0$. Again $D$ is large
in $[Z]$ for any $Z\in \bb$. Take $\beta\colon \ca \bij
(\Q_*^{<\N})^\N$ to be some fixed recursive isomorphism.

Given a $Y=(y_i)\in \bb$ any element
$(\un\lambda_i)\in(\Q_*^{<\N})^\N$ codes a unique infinite
sequence of block vectors (not necessarily consecutive) of $Y$,
which we denote by $(\un\lambda_i)\times Y$. So due to the
minimality of $[e_i]$ there is for any $Y=(y_i)\in \bb$ some
$(\un\lambda_i)\in(\Q_*^{<\N})^\N$ such that
$(e_i)\approx(\un\lambda_i)\times Y$( a standard perturbation
argument shows that the basic sequence $(e_i)$ always embeds as a
sequence of finite rational blocks, though not necessarily
consecutive).

 Set $W=\mgdv Y=(y_i)\in \bb \del (y_{2i}) \in  D \og (e_i)\approx
\beta\circ\alpha(y_{2i})\times (y_{2i+1})\mgdh$, which is then a Borel
subset of $\bb$.

We claim that $W$ is large in $\bb$. For suppose that $Z\in \bb$
is given, take some $V=(v_i)\in D\cap [Z]$ and a
$(\un\lambda_i)\in(\Q_*^{<\N})^\N$ such that
$(e_i)\approx(\un\lambda_i)\times (v_{3i+2})$. Choose
$y_{2i}=v_{3i}$ or $y_{2i}=v_{3i+1}$ such that
$\beta\circ\alpha(y_{2i})=(\un\lambda_i)$ and put
$y_{2i+1}=v_{3i+2}$. Then obviously $Y\leq V\leq Z$ and $Y\in D$.

 So by Gowers' theorem  there is for any $\Delta>0$ a winning
strategy $\tau$ for II for producing blocks in $W_\Delta$ in some
$Y=(y_i)\leq (e_i)$. By choosing $\Delta$ small enough and
modifying $\tau$ a bit we can suppose that the vectors of even
index played by II are in $F_0\cup F_1$. So if $\Delta$ is chosen
small enough, a perturbation argument shows that $\tau$ is in fact
a strategy for playing blocks in $W$.

This shows that if $X$ has a minimal subspace, but does not
contain an isomorphic copy of $c_0$,  there are a basic sequence
$(e_i)$, an element $Y\in \bb$ and a continuous function
$\Theta\colon [Y]\til \overline{span}(e_i)^\N$ (defined by
$\Theta(Z)=\beta\circ\alpha((\tau(Z))_{2i})\times((\tau(Z))_{2i+1})$)
such that for all $Z\in [Y]$ we have
$\overline{span}(\Theta(Z))\subset \overline{span}(Z)$ and
$(e_i)\approx \Theta(Z)$. On the other hand, containing a copy of
$c_0$ is evidently a $\Sigma^1_1$ property, so the disjunction of
the two becomes $\Sigma^1_2$.

\ \ Christian Rosendal\\
Mathematics 253-37\\
California Institute of Technology\\
Pasadena, CA 91125\\
USA\\
 rosendal@caltech.edu, rosendal@ccr.jussieu.fr


\begin{thebibliography}{99}
\bibitem{bla}J. Bagaria, J. Lopez-Abad: Weakly Ramsey sets in
Banach spaces. \emph{Adv. Math., Vol. $\fed {160}$ (2001), no. 2}

\bibitem{bclt}J. Bourgain, P. Casazza, J. Lindenstrauss, L. Tzafriri:
Banach spaces with a unique unconditional basis, up to
permutation. \emph{Mem. Amer. Math. Soc. no. 54 (1985)}

\bibitem{fg}V. Ferenczi, E. M. Galego: Equivalence relations which
are Borel reducible to isomorphism between separable Banach
spaces, {\em preprint}.


\bibitem{fr}V. Ferenczi, C. Rosendal: On the number of non isomorphic
subspaces of a Banach space, \emph{preprint}.

\bibitem{fr2}V. Ferenczi, C. Rosendal: Ergodic Banach spaces,
\emph{preprint}.


\bibitem{g}T. Gowers: An infinite Ramsey theorem and some Banach-space
 dichotomies, \emph{Ann. of Math. 2, Vol. $\fed {156}$,
(2002), no.3}


\bibitem{GM2}T. Gowers, B. Maurey: Banach spaces with small spaces
of operators, \emph{Mathematische Annalen, Vol. $\fed{307}$,
(1997)}

\bibitem{hkl}L. Harrington, A. Kechris, A. Louveau: A Glimm-Effros
dichotomy for Borel equivalence relations,  \emph{Journal of the
Am. Math. Soc., Vol. $\fed{3}$, no. 4 (1990)}

\bibitem{hms}L. Harrington, D. Marker, S. Shelah: Borel orderings.
\emph{Trans. Amer. Math. Soc., Vol. $\fed{310}$, no. 1 (1988)}


\bibitem{HK}G. Hjorth, A. Kechris: Analytic equivalence relations
and Ulm-type classifications, \emph{Journal of symb. logic. Vol.
$\fed {60}$ (1995)}

\bibitem{k}N. Kalton: A remark on Banach spaces isomorphic to their
squares, \emph{Contemp. Math. $\fed {232}$, Amer. Math. Soc.
(1999)}

\bibitem{ke}A. Kechris: Classical Descriptive set theory,
\emph{Springer, New York (1995)}

\bibitem{kj}R. Komorowski, N. Tomczak-Jaegermann: Banach
spaces without local unconditional structure, \emph{Israel J.
Math. $\fed{89}$ (1995), no. 1-3}. Erratum to: "Banach spaces
without local unconditional structure", \emph{Israel J. Math.
$\fed{105}$ (1998)}.

\bibitem{la}J. Lopez-Abad: Coding into Ramsey sets,  \emph{preprint}.

\bibitem{os}E. Odell, T. Schlumprecht: The distortion problem.
\emph{Acta Math., Vol. $\fed{173}$, no. 2}.

\end{thebibliography}
\end{document}